\documentclass[Afour,sagev,times]{sagej}
\usepackage[superscript,biblabel]{cite}
\usepackage{moreverb,url}
\usepackage{amsmath,amssymb,amsfonts}
\usepackage{algorithm}
\usepackage{graphicx} 
\usepackage{textcomp}
\usepackage{multicol, blindtext}
\usepackage{xcolor}
\usepackage{caption}
\usepackage{lipsum}
\usepackage{algorithmic}
\usepackage{amsthm}
\usepackage{subfig}
\usepackage{blindtext, subfig}
\usepackage{dblfloatfix}

\theoremstyle{plain}

\theoremstyle{definition}

\theoremstyle{definition}

\theoremstyle{assumption}

\usepackage[colorlinks,bookmarksopen,bookmarksnumbered,citecolor=red,urlcolor=red]{hyperref}

\newcommand\BibTeX{{\rmfamily B\kern-.05em \textsc{i\kern-.025em b}\kern-.08em
T\kern-.1667em\lower.7ex\hbox{E}\kern-.125emX}}

\def\volumeyear{2016}
\begin{document}

\runninghead{Mondal and Tsourdos}

\title{Cooperative missile guidance design using Distributed Nonlinear Dynamic Inversion}

\author{Sabyasachi Mondal\affilnum{1}, Bhaskar Biswas\affilnum{1}, and Venkatraman Renganathan\affilnum{1}}

\affiliation{\affilnum{1}Cranfield University, UK}

\corrauth{Sabyasachi Mondal, Cranfield University, UK}

\email{sabyasachi.mondal@cranfield.ac.uk}

\begin{abstract}
This paper presents a new cooperative guidance algorithm based on Distributed Nonlinear Dynamic Inversion (DNDI) to demonstrate a coordinated missile attack.   

\end{abstract}


\maketitle
\section{Introduction}

Missiles remain a central element of modern combat systems, and the design of their guidance laws is critical to mission success. A guidance law serves as the decision-making "head" of a missile: it determines how control commands are generated from measurements and estimates so the vehicle steers and intercepts a target. While a single missile can effectively engage a target under favorable conditions, the emergence of advanced air defence systems and increasingly complex mission objectives have highlighted the limitations of single-missile operations.\\\\
To overcome these challenges, cooperative guidance has gained increasing attention. Cooperative guidance refers to the coordination of multiple missiles that share information and adjust their trajectories collectively to achieve objectives such as simultaneous impact, prescribed impact angles, target encirclement, or saturation attacks. Such cooperation not only increases the probability of target neutralization but also complicates enemy countermeasures, since intercepting a group of coordinated missiles is far more difficult than defending against a single one. Furthermore, cooperative engagement enables the attack of larger or multiple targets, which may not be feasible with a single missile.\\\\
Several architectures have been proposed to realize cooperative guidance. Among them, the leader–follower approach is widely adopted due to its simplicity. In this arrangement, a designated leader (or a small subset of leaders) computes high-level references—such as a reference intercept time or desired approach geometry—while followers use local feedback and limited communication with neighbors to track those references and maintain coordination. The leader–follower topology simplifies coordination (fewer variables to negotiate), reduces overall communication load, and can accelerate consensus on mission objectives.\\\\
Despite these advantages, the leader–follower strategy suffers from several limitations:
\begin{itemize}
    \item  If the leader is jammed or destroyed by enemy air defence, the followers lose their reference and the mission is compromised even if they remain intact.
   \item Since followers rely solely on the leader’s commands, any small error in the leader’s computation propagates and amplifies through their own dynamics, leading to significant trajectory deviations.
  \item The lack of direct coordination among followers increases the risk of inter-missile collisions, as each follower is unaware of the precise states of its peers.
  \item Communication delays further degrade performance; if the leader’s command is received late, followers may deviate from the planned trajectory, undermining mission objectives.
\end{itemize}
These drawbacks highlight the need for distributed and decentralized guidance architectures, where each missile integrates information from leaders, neighbors, and its own local states to achieve team-level objectives. Such approaches eliminate over-dependence on a single leader, improve resilience against communication delays and failures, and enable cooperative performance through mutual awareness and coordination among missiles.\\\\
Recent literature reflects significant efforts to address these challenges.
Chen et al. \cite{b1} proposed a leader–follower two-stage consensus law—linearized cooperation followed by PNG—to achieve salvo attack while avoiding singularities, but its reliance on a single leader makes it vulnerable to loss or communication failure.
Yu et al. \cite{b2} designed a fully distributed fixed-time sliding-mode guidance that guarantees convergence of range-to-go and relative velocities under switching topologies, though chattering and careful gain tuning remain issues.
Zhai et al. \cite{b3} used a differential-game formulation to handle maneuvering targets with optimal and adaptive cooperative laws, but solving two-point boundary value problems limits real-time applicability.
Zhang et al. \cite{b4} developed a fully distributed adaptive–optimal strategy using relative distance/velocity consensus and adaptive coupling gains to minimize attack error and control cost without global information, but it assumes bounded target acceleration and requires gain tuning.
Fainkich et al. \cite{b5} and \cite{b7} presented a multi-sliding-surface method that drives neighboring interceptors to reduce time-to-go differences and ensures at least one valid collision course, achieving simultaneous interception without preset impact time but subject to sliding-mode chattering.
Lan et al. \cite{b6} introduced a hybrid coevolutionary guidance law that combines natural evolutionary strategy with PNG to directly optimize time–space consensus, offering high adaptability but heavy computational demands for real-time engagement.
Zhao et al. \cite{b8} proposed a cooperative maneuver-penetration law using line-deviation control to merge attack-time consensus with single-missile spiral maneuvers under overload limits, enhancing penetration but requiring command filtering, observer design, and significant computation.
Tang et al. \cite{b9} formulated a CBF–CLF-based two-stage cooperative guidance that solves a real-time quadratic program for distributed state coordination with nonuniform FOV and input constraints, then switches to PNG for terminal interception, improving robustness at the cost of onboard optimization effort.

Motivated by these challenges, this paper proposes a Distributed Nonlinear Dynamic Inversion (DNDI)-based cooperative guidance framework. Nonlinear Dynamic Inversion (NDI) is regarded as a powerful tool that is useful to design controllers for nonlinear plants. The philosophy behind NDI is to use feedback linearization theory to remove the nonlinearities in the plant. Also, the response of the closed-loop plant is similar to a stable linear system. There are many advantages to using an NDI controller. Some of them are closed-form control expression, easy mechanization, global exponential stability, the inclusion of nonlinear kinematics in the plant inversion, minimization of the need for individual gain tuning or gain scheduling. Dynamic inversion is used to design controllers for many applications. Enns et al. \cite{enns1994dynamic} implemented NDI to design a flight controller. A controller for autonomous landing of UAV was described by Singh et al. \cite{singh2009automatic} using NDI. Padhi et al. \cite{padhi2011neuro} presented Partial Integrated Guidance and Control (PIGC)  to describe reactive obstacle avoidance of UAVs using neuro-adaptive augmented dynamic inversion. A formation flying scheme is proposed by Mondal et al. \cite{mondai2015formation} where the NDI controller is used to track the desired attitude commanded by the leader. Caverly et al. \cite{caverly2016nonlinear} used dynamic inversion to control the attitude of a flexible aircraft. Horn et al. \cite{horn2019non} presented a control design study using dynamic inversion. Lombaerts et al.  \cite{lombaerts2019nonlinear} solved an attitude control problem of a hovering quad tiltrotor eVTOL vehicle using an NDI controller. Extending NDI to a distributed cooperative setting allows each missile to invert its own nonlinear dynamics while simultaneously enforcing coordination objectives with neighboring missiles. This combination enhances robustness, ensures mission reliability under uncertainties, and provides the flexibility required for modern cooperative engagement scenarios.



\section{Preliminaries}
A brief description about the topics required for this work is discussed in this section.
\subsection{Consensus tracking of MASs \cite{mondal2022consensus}} 
Let us consider $N$ nonlinear agents connected by a communication topology. The agents (called followers) need to track the trajectory of a leader, $X_L(t)$ which is connected to a few agents of the networked agents. If the followers' states, i.e., $X_i(t);\ i=1,2,\ldots,N$ achieve the consensus and track the leader's states, i.e., if $X_i(t) \rightarrow X_L(t) $ as $t \rightarrow \infty$, the followers are considered to achieve consensus tracking.

\subsection{Graph Theory}
The communication topology is described using a weighted graph, given by $G= \{V, E\}$. The vertices are given by $V=\{v_1, v_2,\ldots, v_N \}$, which denote the agents. The edges is denoted by the set $E\subseteq V\times V$ which represents the communication among the agents. The elements of weighted adjacency matrix $A=[a_{ij}]\in \Re^{N \times N}$ of $G$ are $a_{ij}>0$ if $(v_i,v_j)\in E$, otherwise $a_{ij}=0$. Since there is no self loop, the adjacency matrix $A$ has zero diagonal elements, i.e., $v_i\in V$, $a_{ii}=0$.
 The degree matrix can be given by $D \in \Re^{N\times N}=diag\{d_1\ d_2\  \ldots d_N\}$, where $d_i=\sum_{j\in N_i}a_{ij}$. The Laplacian matrix is written as $L=D-A$. A graph with the property that $a_{ij}=a_{ji}$ is said to be undirected graph. If any two nodes $v_i, v_j\in V$, there exists a path from $v_i$ to $v_j$, then the graph is called a connected graph. In this paper, we suppose that the topology $G$ of the network is undirected and connected. 

\section{Problem Description}
In this paper, the objective is to design a consensus control for a missile salvo attack on a stationary ground target. To avoid an original leader, We consider one virtual leader and $N$ follower missile. The dynamics of $i^{th}$ missile is given as follows.
\begin{eqnarray}
\dot{V}_i &=& a_{ti} \label{dyn3}\\
\dot{\gamma}_i &=& \frac{a_{ni}}{V_i} \label{dyn4}\\
\dot{x}_i &=& V_i \cos \gamma_i \label{kin1}\\
\dot{z}_i &=& V_i \sin \gamma_i \label{kin2}
\end{eqnarray}
$V_i, \gamma_i, x_i, z_i$ are the velocity, flight-path angle, $x$ and $z$ positions of $i^{th}$ missile respectively. The output $Y_i$ is written as
\begin{equation} \label{yi}
Y_i= \begin{bmatrix}
t_{go_i}\\
\dot{\lambda}_i          
\end{bmatrix} 
\end{equation}
where, $t_{go_i}$ is time-to-go and $\dot{\lambda}_i$ is LOS rate of $i_{th}$ missile. The virtual leader has output  
\begin{equation} \label{y}
Y_L= \begin{bmatrix}
\phi(t,t_{go_0})\\
0          
\end{bmatrix} 
\end{equation}
$\phi(t,t_{go_0})$ is a function of current time and $t_{go_0}=\frac{\sum_{i=1}^N t_{go_i}}{N}$ which we can assume as final time of interception. The objectives of the work can be given as follows.
\begin{itemize}
\item[1.] The follower missiles are required to track a virtual time-to-go profile.
\item[2.] Also, the tracking zero LOS rate assures the interception.
\end{itemize}

\section{Consensus tracking protocol using DNDI controller }
This section will provide an insight into the leader-follower DNDI controller, and it is adopted from the paper \cite{mondal2022consensus}. Let us consider the dynamics for $i^{th}$ agent is given by
\begin{eqnarray} 
\dot{X_i}&=&f(X_i)+g(X_i)U_i \label{non_dyn1}\\
      Y_i&=& h(X_i) \label{non_dyn2}
\end{eqnarray}
The state vector of $i^{th}$ the agent is given by $X_i \in \Re^n$. The control input is $U_i\in \Re^m$. The output of $i^{th}$ the agent is $Y_i\in\Re^p$.
Differentiating Eq. \eqref{non_dyn2} we obtain
\begin{eqnarray}
\dot{Y}_i &=& \left[\frac{\partial h}{\partial X_i} \right] \dot{X}_i \nonumber\\
		  &=& \left[\frac{\partial h}{\partial X_i} \right] \left(f(X_i)+g(X_i)U_i \right) \nonumber\\
		  &=& f_Y(X_i)+g_Y(X_i)U_i 
\end{eqnarray}
where, $f_Y(X_i)=\left[\frac{\partial h}{\partial X_i} \right] f(X_i)$, and $g_Y(X_i)=\left[\frac{\partial h}{\partial X_i} \right] g(X_i)$.
Considering the agent (Eqs. \eqref{non_dyn1}-\eqref{non_dyn2}) and leader output, the consensus tracking error of $i^{th}$ agent (scalar $n=1$) is given by
\begin{equation} \label{error}
{e}_i =\sum_{j\in N_i} a_{ij}\left(y_i-y_j\right)+\beta_i\left(y_i-y_L\right)
\end{equation}
Simplifying Eq. \eqref{error} we get
\begin{equation}
{e}_i= \left(d_i+\beta_i\right)y_i-a_i {Y}-\beta_iy_L
\end{equation}
where $Y=[y_1\ y_2\ \ldots y_N]\in \Re^N$, $y_L$ defines the output of a scalar leader agent, and $\beta_i$ shows if $i^{th}$ agent is connected to the leader. The tracking error is given for the agents with output vector ${Y}_i\in \Re^p; p>1$. 
\begin{equation}
E_i= \left(\bar{d}_i+\bar{\beta}_i\right) Y_i-\bar{a}_i \mathbf{Y}-\bar{\beta}_i Y_L\label{ly3}
\end{equation}
where $Y_L\in\Re^p$, $E_i\in \Re^p$, $\bar{d}_i=(d_i \otimes \mathbf{I}_p) \in \Re^{p \times p}$, $\bar{a}_i=(a_i \otimes \mathbf{I}_p) \in \Re^{p \times pN}$, $\bar{\beta}_i=(\beta_i \otimes \mathbf{I}_p)\in \Re^{p \times p}$, and $\mathbf{Y}=[Y^T_1\ Y^T_2\ \ldots\ Y^T_N]^T \in \Re^{pN}$. We enforce the first-order error dynamics as follows. 
\begin{equation} \label{ly4}
\dot{E}_i+K_i E_i=0
\end{equation}
Differentiation of Eq. (\ref{ly3}) gives
\begin{eqnarray}
\dot{E}_i &=& \left(\bar{d}_i+\bar{\beta}_i\right) \dot{Y}_i - \bar{a}_i \dot{\mathbf{Y}}-\bar{\beta}_i \dot{Y}_L \nonumber\\
          &=& \left(\bar{d}_i+\bar{\beta}_i\right) \big[f_Y(X_i)+g_Y(X_i)U_i \big]-\bar{a}_i \dot{\mathbf{Y}}-\bar{\beta}_i \dot{Y}_L
\end{eqnarray} 
The expressions of $E_i$ and $\dot{E}_i$ are substituted in Eq. (\ref{ly4}) to obtain\\
$\left(\bar{d}_i+\bar{\beta}_i\right)\big[f_Y(X_i)+g_Y(X_i)U_i\big]$
\begin{equation} \label{err_dyn1}
-\bar{a} \dot{\mathbf{Y}}-\bar{\beta}_i \dot{Y}_L + K_i \left(\left(\bar{d}_i+\bar{\beta}_i\right) Y_i-\bar{a}_i \mathbf{Y}-\bar{\beta}_i Y_L\right)=0
\end{equation}
Control $U_{i}$ of $i^{th}$ agent is obtained by simplifying Eq. \eqref{err_dyn1} as follows.\\
$U_{i} = (g_Y(X_i))^{-1} \Big[-f_Y(X_i)+\left(\bar{d}_i+\bar{\beta}_i\right)^{-1} \Big(\bar{a}_i \dot{\mathbf{Y}}+\bar{\beta}_i \dot{Y}_L$
\begin{equation} \label{cont_cons_di}
-K_i \left(\left(\bar{d}_i+\bar{\beta}_i\right) Y_i - \bar{a}_i \mathbf{Y}-\bar{\beta}_i Y_L\Big)\right)\Big]
\end{equation} 

\section{Guidance design using DNDI}
The 2D engagement dynamics between a stationary target and $i^{th}$ interceptor can be given as follows.
\begin{eqnarray}
\dot{r}_i &=& -V_i \cos\phi_i \label{dyn1}\\
\dot{\lambda}_i &=& -\frac{V_i \sin \phi_i}{r_i} \label{dyn2}
\end{eqnarray}
where, $\phi_i = \gamma_i-\lambda_i$. Eqs. \eqref{dyn1}-\eqref{dyn2} are differentiated to obtain 
\begin{eqnarray}
\ddot{r}_i &=& r_i\dot{\lambda}_i^2+u_{i_1} \label{dyn5}\\
\ddot{\lambda} &=& -\frac{2\dot{r}_i \dot{\lambda}_i}{r_i}- \frac{u_{i_2}}{r_i}\label{dyn6}
\end{eqnarray}
Where, $u_{i_1}= \cos\phi_i a_{t_i}-\sin \phi_i a_{n_i}, u_{i_2}=\sin\phi_i a_{t_i}+\cos \phi_i a_{n_i}$, i.e., 
\begin{eqnarray}
\begin{bmatrix}
u_{i_1}\\
u_{i_2}
\end{bmatrix}=
\begin{bmatrix}
\cos\phi_i & -\sin \phi_i\\
\sin\phi_i & \cos \phi_i
\end{bmatrix}
\begin{bmatrix}
a_{t_i}\\
u_{n_i}
\end{bmatrix}
\end{eqnarray}
Let us consider the states as $x_{i_1}=r_i, x_{i_2}=\dot{r}_i, x_{i_3}=\lambda, x_{i_4}=\dot{\lambda}$.  Eqn. \eqref{dyn5} and \eqref{dyn6} are written as 
\begin{eqnarray}
\dot{x}_{i_1} &=& x_{i_2} \label{ss1}\\
\dot{x}_{i_2} &=& x_{i_1}x_{i_4}^2-u_{i_1} \label{ss2}\\
\dot{x}_{i_3} &=& x_{i_4} \label{ss3}\\
\dot{x}_{i_4} &=& -\frac{2x_{i_2} x_{i_4}}{x_{i_1}}- \frac{u_{i_2}}{x_{i_1}} \label{ss4}
\end{eqnarray}
Eqns. \eqref{ss1}-\eqref{ss4} can be written in the form $\dot{X}_i=f(X_i)+g(X_i)U_i$, where, $X_i=[x_{i_1} \ x_{i_2} \ x_{i_3} \ x_{i_4} ]$, $U_i=[u_{i_1} \ u_{_2}]$.
\begin{eqnarray}
f(X_i) &=& \begin{bmatrix}
x_{i_2}\\
x_{i_1}x_{i_4}^2-u_{i_1}\\
x_{i_4}\\
-\frac{2x_{i_2} x_{i_4}}{x_{i_1}}- \frac{u_{i_2}}{x_{i_1}}
\end{bmatrix}
\end{eqnarray}
and 
\begin{eqnarray}
g(X_i) &=& \begin{bmatrix}
0 & 0\\
-1 &0\\
0 &0\\
0 &- \frac{1}{x_{i_1}}
\end{bmatrix}
\end{eqnarray}
The output $Y_i$ is written as
\begin{equation} \label{y}
Y_i= \begin{bmatrix}
t_{go_i}\\
\dot{\lambda}          
\end{bmatrix} 
\end{equation}
$t_{go_i}$ can be expressed as 
\begin{equation} \label{tgo}
t_{go_i} = -\frac{x_{i_1}}{\dot{x}_{i_1}}
\end{equation}
Differentiating Eq. \eqref{tgo} we get
\begin{eqnarray}\label{tgod}
\dot{t}_{go_i}=-1+\frac{x_{i_1}^2x_{i_4}^2}{x_{i_2}^2}-\frac{x_{i_1}}{x_{i_2}^2}u_{i_1}
\end{eqnarray}
The leader's $t_{go}$ is modelled as a straight line and LOS rate is $0$. Therefore, the leader's output $Y_L$ is given by
\begin{equation}
Y_L= \begin{bmatrix}
\delta t+t_{go_0}\\
0          
\end{bmatrix} 
\end{equation}
where, $\delta<0 $ is the slope of the line. 
\section{Results}
\begin{figure}[H]
\centering
\includegraphics[width=3.5in]{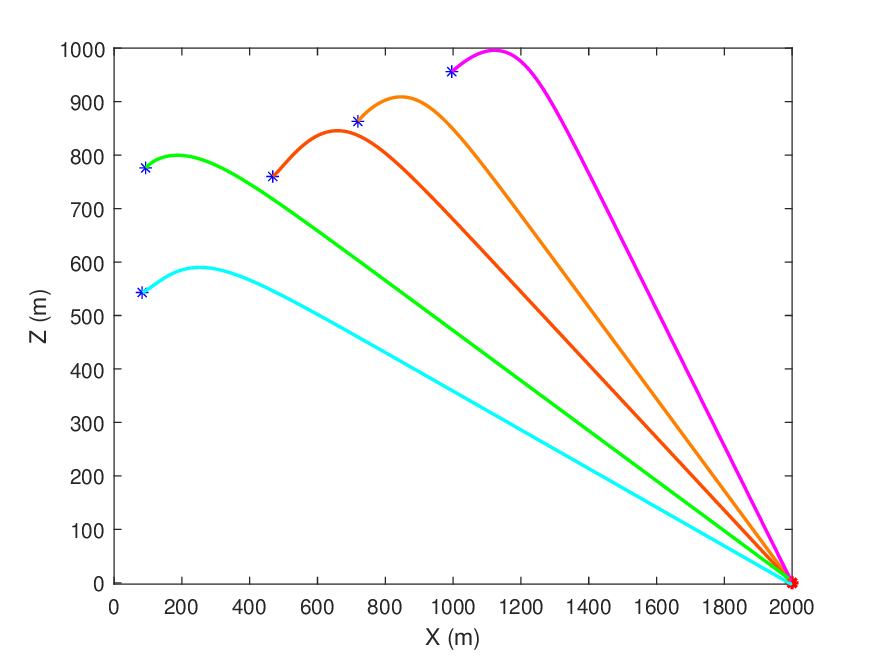}
\caption{Trajectory}
\label{traj}
\end{figure}

\begin{figure}[H]
\centering
\includegraphics[width=3.5in]{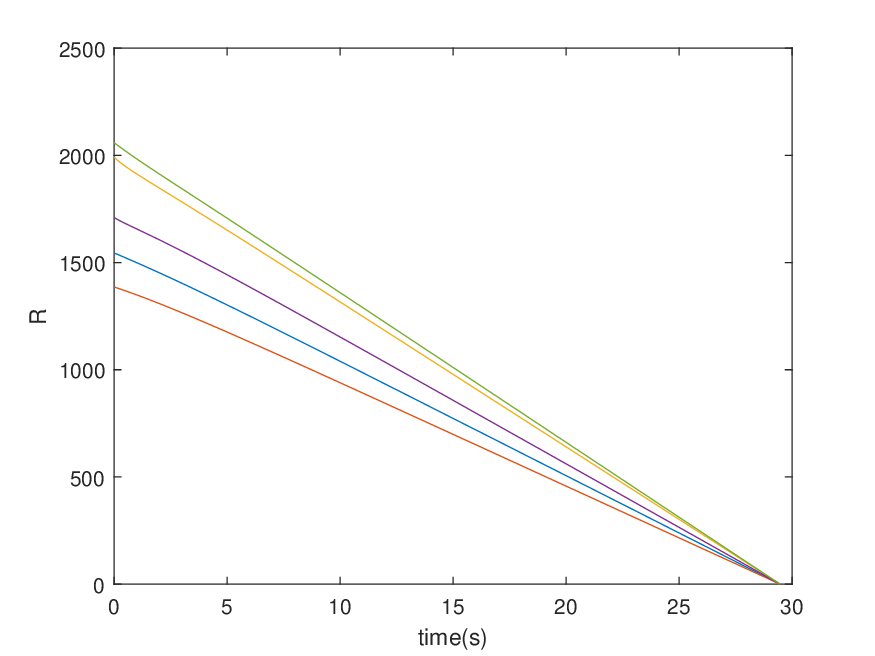}
\caption{Relative distances R.}
\label{r}
\end{figure}

\begin{figure}[H]
\centering
\includegraphics[width=3.5in]{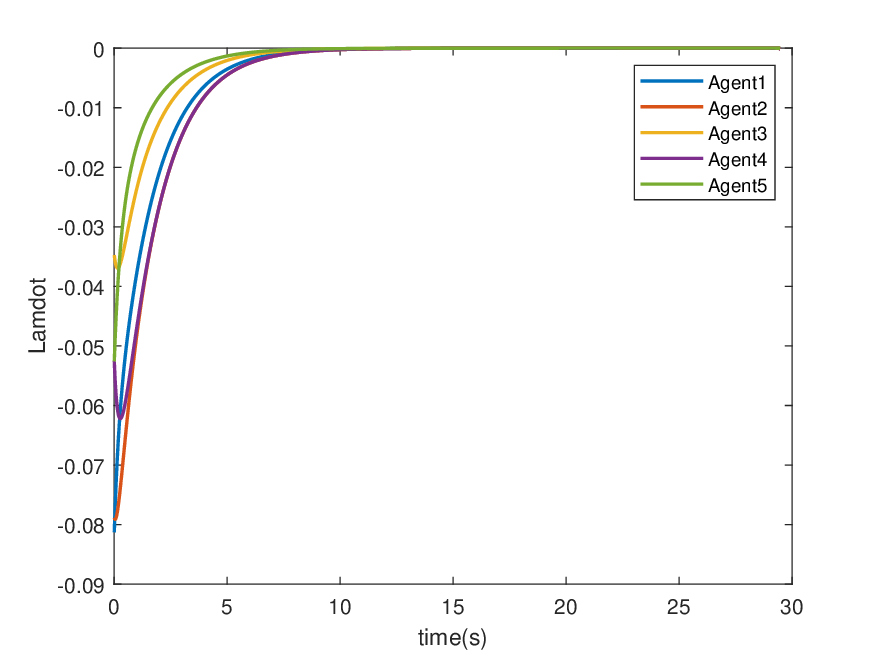}
\caption{LOS rate $\dot{\lambda}$}
\label{lamdot}
\end{figure}

The simulation results are provided here. the target coordinate is fixed at $(0,2000)$. Trajectories of the follower missiles are shown in Fig. \ref{traj}. The target is shown as red dot. Initial positions of the interceptors are shown as blue starts. The missiles hit the target simultaneously which can be verified by the Fig. \ref{tgo} where the missiles achieved the consensus and tracks the leader $t_{go}$ shown in black dotted line. Moreover, the relative distances (Fig. \ref{r}) converge to zero simultaneously. The LOS rate ($\dot{\lambda}$ in Fig. \ref{lamdot}) of missiles achieved zero which show the interception. Therefore, our objectives are achieved using the leader-follower consensus tracking controller designed by DNDI. 
\begin{figure}[H]
\centering
\includegraphics[width=3.5in]{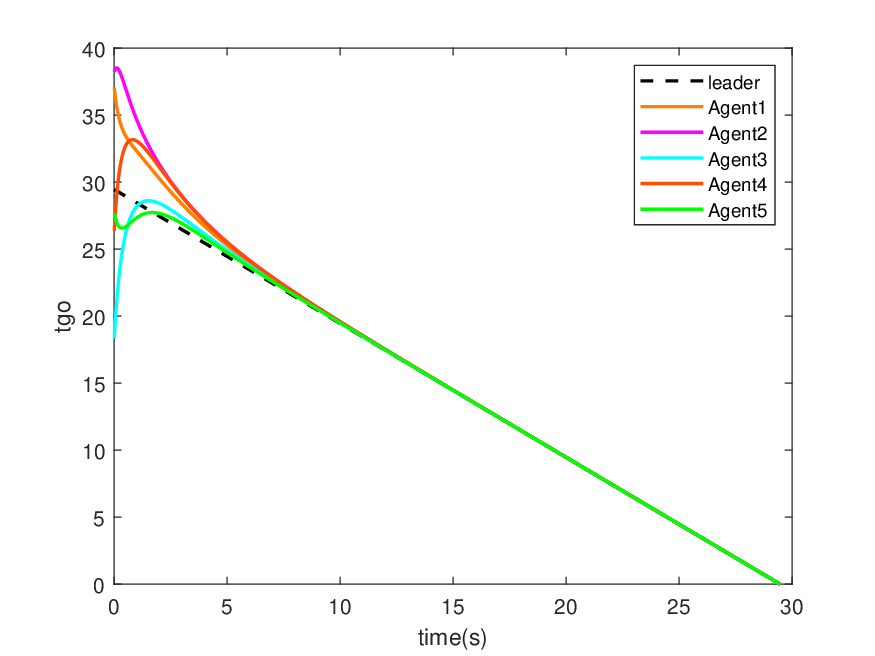}
\caption{Time-to-go}
\label{tgo}
\end{figure}


\section{Conclusion}
A cooperative guidance algorithm has been presented. The effectiveness of the proposed guidance algorithm has been demonstrated using the simulation study. The 

\section*{Acknowledgment}



\bibliographystyle{SageV}
\bibliography{Reference}

\end{document}